\documentclass[11pt,reqno]{amsart}
\usepackage{graphicx, amssymb, color,pdfsync}
\usepackage[all,cmtip]{xy}
\numberwithin{equation}{section}
\usepackage{mathrsfs}
\usepackage{hyperref}

\hypersetup{
    colorlinks=true,       
    linkcolor=blue,          
    citecolor=blue,        
    filecolor=blue,      
    urlcolor=blue           
}
\usepackage{tikz}
\usetikzlibrary{matrix,arrows}

\usepackage[switch]{lineno}
\usepackage{yfonts}
\advance\hoffset by -0.9 truecm

\begin{document}
\newcommand{\s}{\vspace{0.2cm}}

\newtheorem{theo}{Theorem}
\newtheorem{prop}{Proposition}
\newtheorem{coro}{Corollary}
\newtheorem{lemm}{Lemma}
\newtheorem{claim}{Claim}
\newtheorem{example}{Example}
\theoremstyle{remark}
\newtheorem*{rema}{\it Remarks}
\newtheorem*{rema1}{\it Remark}
\newtheorem*{defi}{\it Definition}
\newtheorem*{theo*}{\bf Theorem}
\newtheorem*{coro*}{Corollary}

%
%

\title[Quasiplatonic Riemann surfaces with complex multiplication]{A note on Jacobians of quasiplatonic Riemann \\surfaces with complex multiplication}
\date{}

\author{Sebasti\'an Reyes-Carocca}
\address{Departamento de Matem\'atica y Estad\'istica, Universidad de La Frontera, Avenida Francisco Salazar 01145, Temuco, Chile.}
\email{sebastian.reyes@ufrontera.cl}

\thanks{Partially supported by Fondecyt Grants 11180024, 1190991 and Redes Grant 170071}

\begin{abstract} Let $m \geqslant 6$ be an even integer.  In this short note we prove that  the Jacobian variety of a  quasiplatonic Riemann surface with associated group of automorphisms isomorphic to $C_2^2 \rtimes_2 C_m$ admits complex multiplication. We then extend this result to provide a criterion under which the Jacobian variety of a  quasiplatonic Riemann surface admits complex multiplication.
\end{abstract}
\maketitle

\section{Introduction}
A simple complex polarized abelian variety $A$ of dimension $g$ is said to admit
complex multiplication if its rational endomorphism algebra  $$E=\mbox{End}(A) \otimes_{\mathbb{Z}} \mathbb{Q}$$is a number field of degree $2g$. In this case $E$ is a CM field; namely, a totally imaginary quadratic extension of a totally real field of degree $g$. If $A$ is not simple then, by Poincar\'e Reducibility theorem, there exist pairwise non isogenous simple abelian varieties $A_1,\ldots, A_s$ and positive integers $n_1, \ldots, n_s$ in such a way  that
$$A \sim A_1^{n_1} \times \cdots \times A_s^{n_s}$$where $\sim$ stands for isogeny. By definition, $A$
admits complex multiplication if each simple factor $A_j$ does.

\s

Let $X$ be a compact Riemann surface (or, equivalently, a complex algebraic curve) and let $JX$ denote  its Jacobian variety. Classical examples of compact Riemann surfaces with Jacobian variety admitting complex multiplication are Fermat curves and their quotients. However, in general, it is a difficult task to decide whether or not the Jacobian variety of a given compact Riemann surface admits complex multiplication, and much less is known about their distribution in the moduli space $\mathcal{A}_g$ of principally polarized abelian varieties; see, for example, \cite{OortR}. As a matter of fact, a well-known  conjecture due to Coleman predicts that if $g \geqslant 4$ then the number of isomorphism classes of compact Riemann surfaces of genus $g$ with Jacobian variety admitting complex multiplication is finite. In spite of the fact that this conjecture has been proved to be false for $g \leqslant 7,$ currently it still remains  as an open problem for $g \geqslant 8.$ This conjecture is closely related to important open problems of Shimura varieties, special points in the Torelli locus of $\mathcal{A}_g$ and the theory of unlikely intersections.

\s

If the Jacobian variety of a compact Riemann surface $X$ admits complex multiplication then $X$ can be defined, as complex algebraic curve, over a number field; see \cite{ST}. In part due to this fact, there has been an increase in the interest of these compact  Riemann surfaces, particularly in their applications to number theory and arithmetic geometry.

\s

By the classical Belyi's theorem \cite{Belyi}, a compact Riemann surface $X$ can be defined over a number field if and only if there exists a holomorphic map $$\beta : X \to \mathbb{P}^1$$ with at most three critical values; the pair $(X, \beta)$ is called a  Belyi pair. Possibly the more interesting examples of Belyi pairs are the  regular ones: namely, those for which $\beta$ is given by the action of a group of automorphisms of $X$. In this case, $X$ is known to be  quasiplatonic (or to have many automorphisms); that is, it cannot be deformed non-trivially in the moduli space together with its automorphism group.  

\s

Oort in \cite[p.18]{Oort} considered quasiplatonic Riemann surfaces and discussed the problem of determining which among them have Jacobian variety admitting complex multiplication. For genus at most four, this problem was completely solved by Wolfart in \cite{JW}. Later, M\"{u}ller and Pink in \cite{MP} and Obus and Shaska in \cite{Shaska} considered the hyperelliptic and superelliptic situations respectively, and succeeded in determining which among them have Jacobian variety admitting complex multiplication.

\s

A different approach can be done by considering regular Belyi pairs $$\beta: X \to \mathbb{P}^1 \cong X/H$$ whose covering groups $H$ share a common property or have the same algebraic structure. In this direction, it was proved in   \cite{JW} that if $H$ is abelian then $JX$ admits complex multiplication, and in \cite{clr} the same conclusion was obtained for two infinite series of compact Riemann surfaces arising as quotients of regular Belyi curves with a metacyclic group of automorphisms.

\s

In this short note we consider a infinite series of quasiplatonic Riemann surfaces with associated covering group isomorphic to the  semidirect product $$G_m:=\langle a, b, t: a^2=b^2=(ab)^2=t^m=1, tat^{-1}=a, tbt^{-1}=ab \rangle \cong C_2^2 \rtimes_2 C_m$$where $m \geqslant 6$ is even integer. Based on the classification and description obtained in \cite{nos} for quasiplatonic Riemann surfaces with action of semidirect products of the form  $C_2^2 \rtimes C_m$, we prove the following result.

\begin{theo*}  Let $m \geqslant 6$ be an even integer. If $(X, \beta)$ is a regular Belyi pair with associated covering group isomorphic to $G_m$ then the Jacobian variety $JX$ admits complex multiplication. 
\end{theo*}

The proof of the theorem --which is rather simple and based on the classical theory of covering of Riemann surfaces-- is done in Section \S\ref{dos}. Then, in Section \S\ref{tres} we extend the arguments used to prove the theorem to provide a criterion under which the Jacobian variety of a quasiplatonic Riemann surface admits complex multiplication. Finally, we end this short note by recalling a couple of observations in Section \S\ref{s4}.

\section{Proof of the theorem}\label{dos}
Let $(X, \beta)$ be a regular Belyi pair with associated covering group isomorphic to $G_m$ where $m \geqslant 6$ is even.

\s

{\bf Case A.} Assume $m\equiv 2 \mbox{ mod }4.$ Following \cite[Theorem 1(1)]{nos}, the regular covering map $$X \to X/G_m \cong \mathbb{P}^1$$ramifies over three values marked with  $2,m$ and $2m,$ and the genus of $X$ is $$g_X=m-2.$$ Moreover, as observed in \cite[Subsection \S 4.3.1]{nos}, its Jacobian variety decomposes isogenously as the product  \begin{equation}\label{iso}JX \sim JY^2\end{equation}where $Y=X/\langle a \rangle$ is the quotient Riemann surface represented by the affine algebraic curve $$y^2 = x^{m}-1.$$

Note that $\langle a \rangle$ is a normal subgroup of $G_m$ and the quotient $H=G_m/\langle a \rangle$ is an abelian group of order $2m$ acting as a group of automorphisms of $Y.$ Clearly, the corresponding orbit space $Y/H$ has genus zero and the associated regular covering map $$\beta_{H}: Y \to Y/H$$ramifies over at most three values.  It follows that $(Y, \beta_H)$ is a regular Belyi pair with abelian covering group; thus, by \cite[Theorem 4]{JW}, we obtain that $JY$ admits complex multiplication. The result follows from the isogeny  \eqref{iso}.

\s

{\bf Case B.} Assume $m\equiv 0 \mbox{ mod }4$.  Following \cite[Theorem 1(2b)]{nos}, the regular covering map $$X \to X/G_m \cong \mathbb{P}^1$$ramifies over three values marked with $2, m$ and $m$  and the genus of $X$ is $$g_X=m-3.$$ Moreover, as observed in \cite[Subsection \S 4.3.3]{nos}, its Jacobian variety decomposes isogenously as the product\begin{equation}\label{iso3}JX \sim JY \times JZ^2\end{equation}where $Y=X/\langle a \rangle$ and $Z=X/ \langle b \rangle$ are the quotient Riemann surfaces represented by the affine algebraic curves $$y^2 = x^{m}-1 \,\, \mbox{ and }\,\, y^2=x^{\frac{m}{2}}-1$$respectively; their genera  are $g_Y=\tfrac{m}{2}-1$ and $g_Z=\tfrac{m}{4}-1.$

\s

We argue analogously as done in the case {\bf A} to ensure that $JY$ admits complex multiplication. Besides, in order to prove that $JZ$ also does, define $$\iota(x,y)=(x,-y) \,\, \mbox{ and } \,\,\,\tau(x,y)=(\mbox{exp}(\tfrac{4 \pi i}{m})x,y)$$ and let $K=\langle \iota, \tau \rangle \cong C_2 \times C_{m/2}.$ We observe that the abelian group $K$ satisfies $$|K| > 4(g_Z-1) \,\, \mbox{for all }m \geqslant 8.$$Then, by the classification of large abelian groups of automorphisms of compact Riemann surfaces given in \cite[Theorem 3.1]{Lomuto}, we see that the branched regular covering map  $$\beta_K: Z \to Z/K \cong \mathbb{P}^1$$ramifies over three values, marked with 2, $\tfrac{m}{2}$ and $\tfrac{m}{2}$.   Thus, again by \cite[Theorem 4]{JW}, we conclude that $JZ$ admits complex multiplication and the result follows from the isogeny  \eqref{iso3}.

\section{A generalization} \label{tres}

 Let $(X, \beta)$ be a regular Belyi pair and let $G$ denote the associated covering group. Consider a collection $$\{H_1, \ldots, H_s\}$$  of proper non-trivial subgroups of $G.$ Let $Y_i$ denote the quotient Riemann surface $X/H_i$ and let $g_i \neq 0$ denote its genus, for each $i \in \{1, \ldots, s\}.$

\s

Assume the existence of  positive integers $n, n_1, \ldots, n_s$ in such a way that $$JX^n \sim JY_1^{n_1} \times \cdots \times JY_s^{n_s}$$(we point out that conditions under which an isogeny as above can be obtained were determined, for example, in \cite{KR} and later generalized in  \cite{KRnos}).

\s

Consider the following statements:

\begin{enumerate}
\item[{\bf A.}] $H_i$ is a normal subgroup of $G$ and  $G/H_i$ is abelian.
\item[{\bf B.}] $Y_i$ admits a large abelian group $K_i$ of automorphisms (namely, its order is strictly  greater than $4(g_i-1)$) with only one exception: $K_i \cong C_6$ and $$Y_i \to Y_i/K_i$$ramifies over four values; two marked with 2 and two marked with 3.
\end{enumerate}

The arguments employed in the proof of the theorem are naturally generalized to provide the following criterion. With the same notations:

\s

{\bf Criterion.} If for each $i \in \{1, \ldots, s\}$ either $H_i$ satisfies the statement ${\bf A}$ or $Y_i$ satisfies statement ${\bf B}$, then $JX$ admits complex multiplication.

\s

It is worth to mention that the statement {\bf B} can be restated in a weaker manner. Indeed, the same conclusion is obtained if we ask $Y_i$ to be endowed with a quasi-large abelian group $K_i$ of automorphisms (namely, its order is strictly  greater than $2(g_i-1)$) not belonging to one of the $22$ exceptional cases listed in \cite[Table 2]{quasi}.  
\section{Remarks} \label{s4}

\subsection*{Remark 1} We should mention that the criterion is, as expected, rather restrictive. However, it provides a different approach to find new examples of Jacobian varieties admitting complex multiplication. In addition, it is worth recalling that a shorter proof of the theorem can be obtained by noticing that  $Y$ and $Z$ in the theorem are hyperelliptic. Nevertheless, as our proof is based on significantly simpler arguments which do not depend on the hyperellipticity of the involved quotients, its generalization could be used for a possibly wider range of cases.

\subsection*{Remark 2}  In  \cite{St} Streit provided a representation theoretic sufficient condition for the Jacobian variety of a quasiplatonic Riemann surface to admit complex multiplication. Concretely, with the previous notations, if $\mbox{S}^2(\rho_a)$ denotes the symmetric square representation of the analytic representation $\rho_a$ of $G$ and $1$ stands for the trivial representation of $G$ then $$\langle \mbox{S}^2(\rho_a), 1 \rangle_G=0 \, \implies \, JX \mbox{ admits complex multiplication.}$$

After routine computations, one sees that the previous criterion allows to conclude that $JX$ admits complex multiplication provided that $n \equiv 2 \mbox{ mod } 4.$ However, this criterion does not provide conclusion if $n \equiv 0 \mbox{ mod } 4.$

\subsection*{Acknowledgements} The author is grateful to the referee for suggesting useful improvements to the article, and to Jennifer Paulhus and Anita M. Rojas for valuable conversations and for sharing their computer routines with him.

\end{document}